\documentclass{article}
\usepackage{graphicx} % Required for inserting images
\usepackage{amsmath}

\usepackage{geometry}
\usepackage{tabularx}
\usepackage{enumitem}
\usepackage{geometry}
\usepackage{array}
\usepackage{multirow}
\usepackage[toc]{appendix}
\usepackage{xpatch}

\title{Comparative Study of Sand Drawings in Oceania and Africa}
\author{Linbin Wang, Rowena Ball, Hongzhang Xu}

\usepackage{caption}

\usepackage{animate}
\usepackage{graphicx}

\graphicspath{ {figures/} }
\usepackage{array}

\usepackage{biblatex}

\DeclareNameAlias{sortname}{last-first}

\renewbibmacro{in:}{} % Remove "In:" before journal name

\renewbibmacro*{volume+number+eid}{%
  \printfield{volume}%
  \setunit{\addspace\bibstring{number}\addnbspace}%
  \printfield{number}%
  \setunit{\addcomma\space}%
  \printfield{eid}}

\DeclareFieldFormat[article]{volume}{Volume~\mkbibbold{#1}}

\xpatchbibmacro{date+extrayear}{%
  \printtext[parens]%
}{%
  \printtext%
}{}{}

\DefineBibliographyStrings{english}{%
  page = {page},
  pages = {pages},
  volume = {Volume},
  number = {Number},
}
  
\begin{document}
\widowpenalty10000
\clubpenalty10000

\date{}
\maketitle
\indent
\pdfoutput=1

\begin{abstract}
People typically consider only European mathematics as orthodox, often intentionally or unintentionally overlooking the existence of mathematics from non-European societies. Inspired by Maria Ascher’s two well-known papers on sand drawings in Oceania and Africa, this paper focuses on the strong link between modern mathematics and the mathematics behind the sand drawings. Beginning with a comparison of the geography, history, and the cultural context of sand drawings in Oceania and Africa, we will examine shared geometric features of European graph theory and Indigenous sand drawings, including continuity, cyclicity, and symmetry. The paper will also delve into the origin of graph theory, exploring whether the famous European mathematician Leonhard Euler, who published his solution to the Königsberg bridge problem in 1736, was the true inventor of graph theory. The potential for incorporating sand drawings into the school curriculum is highlighted at the end. Overall, this paper aims to make readers realise the importance of ethnomathematics studies and appreciate the intelligence of Indigenous people.   
\end{abstract}
\hspace{10pt}
% Keywords command
\providecommand{\keywords}[1]
{
  \small	
  \textbf{\textit{Keywords---}} #1
}


\section{Introduction}

In the realm of mathematics, the contributions of colonised peoples were often ignored or devalued \cite{Josepha20}. It wasn't until the 1970s and 1980s that a growing resistance emerged among teachers and mathematics educators in developing countries, and later in other countries as well, against the racist and (neo)colonial prejudices related to mathematics and against the Eurocentrism in mathematics and its history. It was stressed that beyond the “imported school mathematics," there are, and continue to be, other forms of mathematics \cite{Gerdesa29}.

This is often attributed to the dominance of Europe and its cultural dependencies \cite{Josepha20} — countries, notably the United States, Canada, Australia, and New Zealand, which are inhabited mainly by populations of European origin or which have historical and cultural roots similar to those of European peoples — in global affairs over the last four hundred years. This dominance is too often reflected in the character of some of the historical writing produced by Europeans in the past. Where other people appear, they do so in a transitory fashion whenever Europe has chanced in their direction. Consequently, the history of Africans and other Indigenous peoples worldwide often seems to commence only after their encounter with Europe, resulting in the oversight of their pre-colonial contributions to scientific fields such as mathematics.

With the increasing recognition of ethnomathematics studies, various concepts have been proposed by distinguished ethnomathematicians to contrast with “academic mathematics" or “school mathematics" (i.e. the curriculum transplanted from elsewhere) : \textit{indigenous mathematics} (Gay and Cole; Lancy); \textit{sociomathematics} (Zaslavsky); \textit{informal mathematics} (Posner); \textit{spontaneous mathematics} (D'Ambrosio); \textit{oral mathematics} (Canaher, Carraher and Schliemann); \textit{oppressed mathematics} (Gerdes); \textit{hidden or frozen mathematics} (Gerdes); \textit{folk mathematics} (MellinOlsen); \textit{non-standard mathematics} (Carraher et al.; Gerdes; Harris); and \textit{mathematics codified in know-how} (Ferreira) \cite{Gerdesa29}. 

According to Vandendriessche and Silva \cite{vandendriesschea2}, among all the outstanding ethnomathematicians in history,  the British anthropologist John Layard (1891-1974) and the young Cambridge student Arthur Bernard Deacon (1903-1927) were the first ethnographers to document the sand drawing practice (the practice of local people drawing geometrical figures using one finger on the sand to create continuous patterns) in the New Hebrides (which became the Republic of Vanuatu in 1980). Under the leadership of the Cambridge anthropologist Alfred Cort Haddon (1855-1940), Deacon carried out ethnographic research between 1926 and 1927, mainly in the two Seniang and Mewun societies on the southwest coast of the island of Malekula. He collected a large number of these sand drawings and recorded each drawing precisely, giving all the information allowing its reproduction. In the 1980s, American mathematician Marcia Ascher (1935-2013), one of the founding members of ethnomathematics, chose to analyse the corpus of sand drawings collected with precision by Layard and Deacon \cite{vandendriesschea2}.

This article is motivated by Marcia Ascher’s two well-known papers, “Graphs in Cultures: A Study in Ethnomathematics" \cite{ascher1987a3} and “Graphs in Cultures (II): A Study in Ethnomathematics \cite{ascher1988a4}." These two papers highlight the presence of mathematics in various Indigenous cultural contexts, suggesting that mathematics is a natural practice carried out instinctively by all humans.

\section{Methods}

The paper aims to compare and contrast “Graphs in Cultures: A Study in Ethnomathematics” and “Graphs in Cultures (II): A Study in Ethnomathematics," showcasing the differences and correlations between sand drawings discovered in New Ireland and Malekula in Oceania and those found in Africa among the Bushoong people in Zaire and the Tshokwe people in Angola. 

We first examined the geography, history, and cultural context of each region. Subsequently, we focused on shared attributes between European graph theory and Indigenous sand drawings, including continuity, cyclicity, and symmetry. Afterward, we explored whether the renowned European mathematician Leonhard  Euler, who proposed a solution to the Königsberg bridge problem in 1736, was the pioneer of graph theory and drew inspiration from any Indigenous culture to solve the problem. Additionally, we underscored the potential to integrate sand drawings into the school curriculum. A table was employed to succinctly and efficiently illustrate the comparison (See Table 1).

\begin{table}
  \centering
  \begin{tabularx}{\linewidth}{|l|X|X|}
    \hline
    \textbf{} & \textbf{Study I} & \textbf{Study II} \\
    \hline
    \textbf{Geography} & In Oceania, New Ireland in Papua New Guinea and Malekula in Vanuatu. & In Africa, among the Bushoong culture in Zaire and the Tshokwe culture in northeastern Angola. \\
    \hline
    \textbf{Cultural Meaning} & \begin{itemize}[label=\textbullet, left=0pt, nosep, after=\strut]
                                   \item In Malekula, sand drawings, named after local flora and fauna, are connected to myths and rituals.
                               \end{itemize}
                             & \begin{itemize}[label=\textbullet, left=0pt, nosep, after=\strut]
                                   \item The Bushoong view a sand drawing design as decomposed into different elementary designs, and the name given to a figure is the name associated with its most significant constituent.
                                   \item For the Tshokwe people, sona traditionally depicted proverbs, fables, riddles, and animals, serving as a vital means of knowledge transmission.
                               \end{itemize}\\

    \hline
    \textbf{History} & Sand drawing tradition documented since 1926 when Deacon carried out ethnographic research in Malekula. & Some circumstantial evidence proves that the sand drawing tradition must be ancient and certainly pre-colonial in sub-Saharan Africa. \\
    \hline
    \textbf{Shared Geometrical Features} & \multicolumn{2}{c|}{Continuity, Cyclicity, Symmetry} \\
    \hline
    \textbf{Educational Potential} & In 2019, Vandendriessche and Worwor integrated local activities into math curriculum in North Ambrym, Vanuatu. & Gerdes explored using Tshokwe sand drawings in math education, covering arithmetical relationships, symmetry, geometric algorithms, and generalization.\\
    \hline
  \end{tabularx}
  \caption{Comparison of Studies on Sand Drawings}
  \label{tab:sand_drawings_horizontal}
\end{table}

\section{Results}

\subsection{Geography, History, and Cultural Context}
“Graphs in Cultures: A Study in Ethnomathematics” concentrates on the sand drawings found in New Ireland and Malekula. 

New Ireland is the largest island of New Ireland Province, situated in the north easternmost part of Papua New Guinea. Although there are limited reports on sand drawings in New Ireland, mainly from F. L. S. Bell, these drawings will be discussed to supplement the argument for the universal nature of graph theory in human culture.

Malekula is an island roughly 80 miles long, and lies in the centre of the Northern New Hebrides \cite{layard1972a5}. It is the second largest island of the New Hebrides group and is inhabited by a people speaking dialects of a Melanesian language. Malekula is divided into a number of districts which differ from each other both culturally and linguistically, those which have been studied most closely are Seniang and Mewun in the south-west, and Lambumbu and Lagalag in the narrow part of the island towards the north (see Figure \ref{map of Malekula}) \cite{deacon1934a6}. Exploitation of this island by Europeans began in about 1840. The social and psychological effects of colonial rule and the introduction of European diseases reduced the population by about 55\% between the 1890s and 1930s \cite{gathercole1984a7}.

\begin{figure}
    \centering
    \includegraphics[width=1\linewidth]{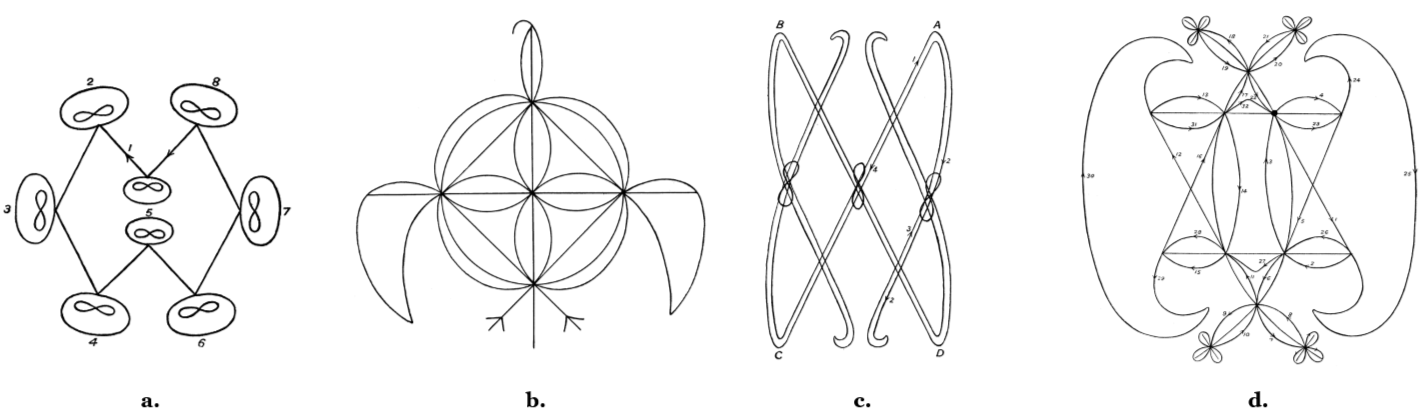}
    \caption{(a) \textit{Nimbip}: A variety of tuber. The figure represents two of these tubers planted in the centre, and on the circumference, six tubers which have sprouted from them. (b) \textit{Nooimb}: A kind of bird. (c) \textit{Dimbuk Temes}: This is a flower head-dress worn in the funeral ritual of the Nalawan ceremonies. (d) \textit{Nevet Tambat}: “The Stone of Ambat." It likely alludes to a story involving the Ambat brothers, important mythological culture heroes of Seniang, in which the mythical personage Nevinbumbaau attempts to kill and eat them. }
    \label{fig.21}
\end{figure}

As mentioned in the introduction, a large part of the graphs used in Ascher’s ethnomathematics study in Malekula are adapted from the field-notes and letters of A. Bernard Deacon, who, during 1926 and 1927, spent a year making an ethnographic study of the natives in South West Bay, Malekula \cite{ascher1987a3}. Deacon's collection reveals that many of the figures are named after flora and fauna in the local environment, while several are related to important myths and rituals. Figure \ref{fig.21} exemplifies each type. In addition to illustrating myths, there are even myths in which the drawing concept plays a significant role. A Malekula tale which incorporates the concept of drawing geometrical figures will be covered later. The art of drawing, called nitus, is handed down from generation to generation, and it requires considerable skill and practice. It is noteworthy that knowledge of the art is entirely limited to men; women, of course, may see the designs \cite{deacon1934a6}.

\begin{figure}
    \centering
    \includegraphics[width=1\linewidth]{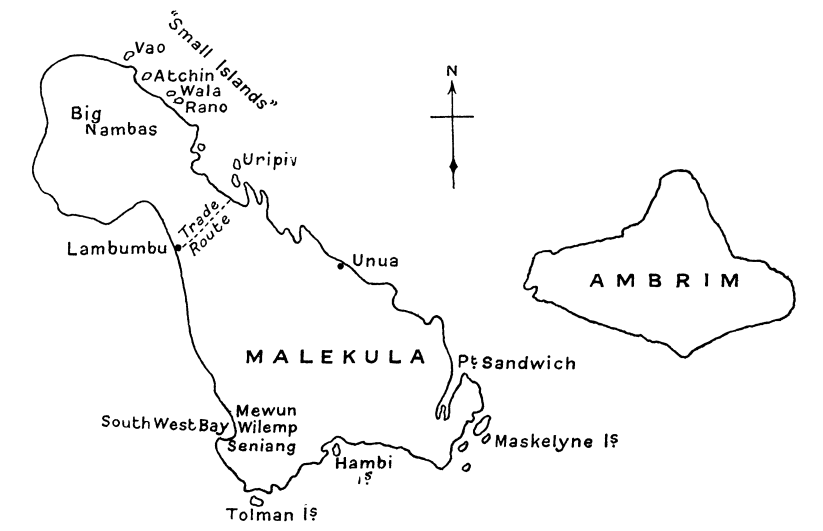}
    \caption{Sketch map of Malekula.}
    \label{map of Malekula}
\end{figure}

In about 1905, a European ethnologist studying among the Bushoong in Africa was challenged by some children to trace figures in the sand. The challenge specified that each line should be traced once and only once without lifting his finger from the ground. At that time, he was unaware that the Bushoong shared concerns with ideas now known as graph theory. The Bushoong are one of the subgroups of the Kuba chiefdom, which became Zaire in 1960 after the colonization by the Belgian government. The Bushoong view a sand drawing design as decomposed into different elementary designs, and the name given to a figure is associated with its most significant constituent \cite{ascher1988a4}. 

On the other hand, in the northeastern region of Angola in Africa, Gerdes \cite{gerdes1988a9} discovered that the sand drawing tradition (called lusona, sona for plural) also exists and belongs to the heritage of the Tshokwe, Lunda, Lwena, Xinge, and Minungo peoples. Among these Southern African peoples, the Tshokwe are especially well-known for their beautiful decorative art, ranging from the ornamentation of plaited mats and baskets, ironwork, ceramics, engraved calabash fruits, and tattoos to paintings on house walls and sand drawings.

When the Tshokwe gathered at their central village places or hunting camps, they usually sat around a fire or in the shade of leafy trees, spending their time in conversations illustrated by drawings in the sand. These drawings depicted proverbs, fables, riddles, animals, and more, playing an important role in transmitting knowledge and wisdom from one generation to the next \cite{gerdes1988a9}. This can be compared to the modern teaching method in class, where teachers use chalk to write down words and sometimes draw graphs on the blackboard to convey knowledge to students. Much like the sand drawing tradition in Malekula, the meaning and execution of simple sona were exclusively reserved for Tshokwe men, often acquired during the intensive schooling phase of circumcision and initiation rites. The meaning and execution of more difficult sona were known only by specialists, the akwa kuta sona (those who know how to draw), who again transmitted their knowledge to their sons \cite{gerdes1988a9}.

\begin{figure}
    \centering
    \includegraphics[width=1\linewidth]{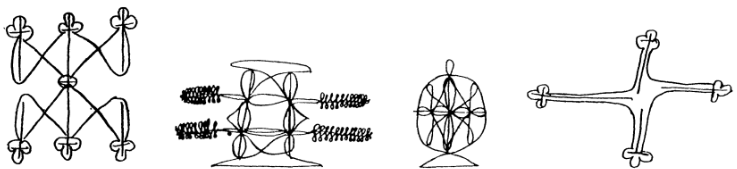}
    \caption{ Sand drawings attached to tales collected from New Hebrides islands. From left to right are: bwat gaviga burabura, ruerue, avua, and tum revrevo. }
    \label{Sand drawings attached to tales}
\end{figure}

Despite the different geographic location, both cultures use sand drawings for cultural expression, storytelling, and knowledge transmission. These stories, passed down from generation to generation, carry the unique history and beliefs of each culture. The narratives exhibit variability. Some tales collected from natives of the islands of New Hebrides have attached picture-drawings (see Figure \ref{Sand drawings attached to tales} \cite{firth1930a10}). Some stories are considered fiction and lack real moral value. For example, the tale involving a boy being pulled down from a tree is regarded as a humorous anecdote \cite{firth1930a10}. 

\begin{figure}
    \centering
    \includegraphics[width=0.22\linewidth]{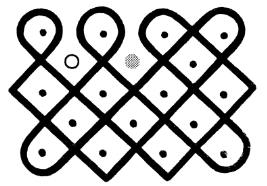}
    \caption{In this drawing, the gray circular point represents the hunter, and the open circular point the dog.}
    \label{hunter and dog}
\end{figure}

Nevertheless, many stories carry a deep philosophical meaning that reflects the wisdom of Indigenous people. A famous Malekula tale \cite{deacon1934a6} symbolizes the cycle of life, death, and potential rebirth, focusing on the renowned warrior Airong who died in Benaur. His ghost attempted to traverse a road near the rock where Temes Savsap guarded a geometrical figure. As Airong had never learned this figure, he sought in vain to pass the rock. When Temes was about to devour him, Airong narrowly escaped. After miraculously returning to life, Airong requested his bow and arrow and defeated Temes. The incorporation of a geometrical figure in the story indicates how mathematical ideas are closely embedded in the local culture, and the way Temes eventually overcame the challenge highlights the significance of courage when facing difficulties in life. Another traditional Tshokwe fable that also conveys a message is “The Hunter and the Dog” \cite{gerdes1988a9}, which underscores the principle of fairness in resource-sharing. Initially, hunter Tshipinda divided the meat with Calala, the dog's owner, leaving the dog with only bones. When Tshipinda sought the dog's help again, the dog refused and advised the hunter to take Calala instead. This thus reinforces the importance of prior agreements and equitable divisions. A corresponding lusona is shown in Figure \ref{hunter and dog}.

\begin{figure}
    \centering
    \includegraphics[width=1\linewidth]{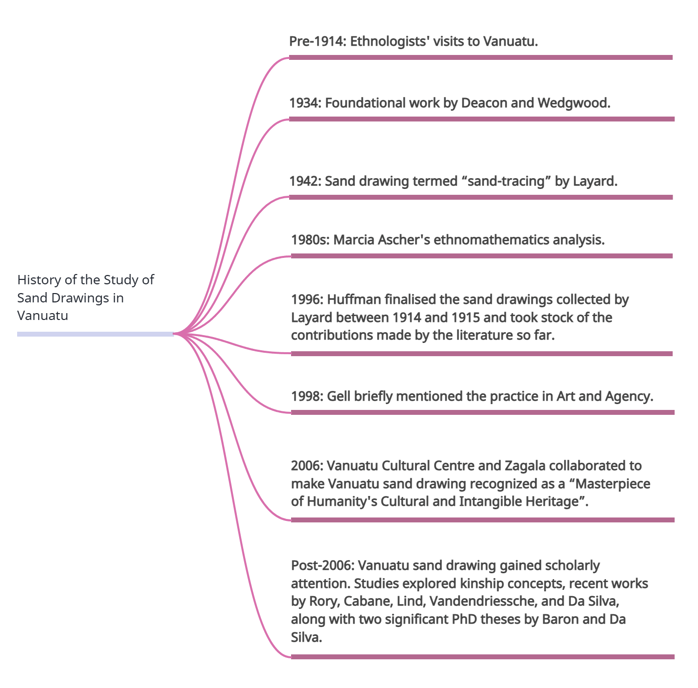}
    \caption{History of the ethnomathematical study of sand drawings in Vanuatu}
    \label{history1}
\end{figure}

\begin{figure}
    \centering
    \includegraphics[width=1\linewidth]{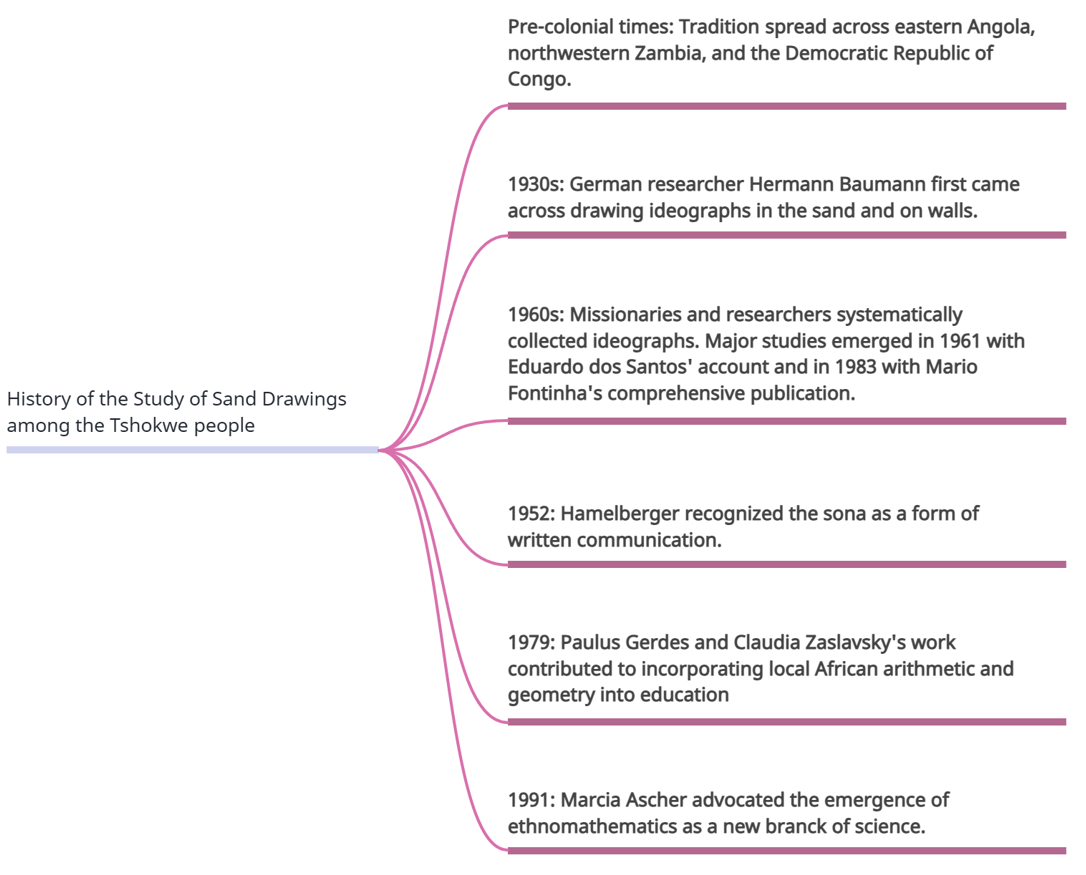}
    \caption{History of the ethnomathematical study of sand drawings in Angola}
    \label{history2}
\end{figure}

Furthermore, the history of the ethnomathematical exploration of sand drawings, both in Vanuatu \cite{Devyldera36} and among the Tshokwe people in Angola \cite{Kubikea30}, is encapsulated in the two diagrams provided below (See Figure \ref{history1} and Figure \ref{history2}). It is evident from these diagrams that the discovery of sand drawings in both regions follows a similar story, showcasing how the growing awareness of this tradition has developed through the ongoing efforts of ethnomathematicians.

\subsection{Geometrical Ideas}

Continuity, cyclicity, and symmetry occur as not merely prevalent drawing themes in the sand drawings found in Oceania and Africa but also serve as foundational concepts in graph theory. A closer examination of the significance of these elements within Indigenous communities provides valuable insights into how Indigenous people perceive and engage with mathematics, contributing to a better understanding of the intersection between traditional practices and mathematical concepts within these communities.

\subsubsection{Continuity}

\begin{figure}
    \centering
    \includegraphics[width=0.5\linewidth]{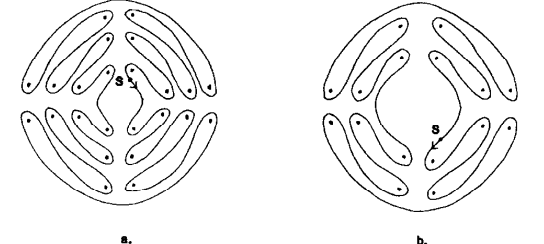}
    \caption{Simple closed curves: (a) Ambat, trying to sharpen a stone on the seashore, being pushed back by the incoming tide. (b) Associated with a secret society. }
    \label{simple closed curves}
\end{figure}

In a letter dated October 1926 from Lambumbu, reporting progress to Dr. Haddon, Deacon wrote: “The whole point of the art is to execute the designs perfectly, smoothly, and continuously; to halt in the middle is regarded as an imperfection” \cite{deacon1934a6}. The idea of drawing continuously is further supported by Ascher's paper “Graphs in Cultures: A Study in Ethnomathematics”, which notes, “Of those for which it is possible to use a single continuous line and the courses are known, 94\% actually were traced that way”. Two examples are listed in Figure \ref{simple closed curves} \cite{ascher1987a3}. Notice that the pair is isomorphic, which is another significant idea in modern mathematics.

\begin{figure}
    \centering
    \includegraphics[width=0.3\linewidth]{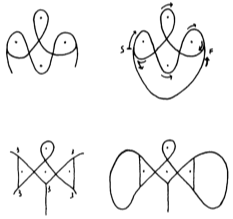}
    \caption{ The upper two graphs represent a “dyahotwa" bird but are not traceable. The graphs at the bottom are an “extension" of the ones above in order to achieve that.
 }
    \label{mono-linear graphs}
\end{figure}

In a similar way, the designs of the Tshokwe sand drawings must be executed smoothly and continuously. Any hesitation or pause on the part of the drawer is interpreted by the audience as imperfection and a lack of knowledge, often met with an ironic smile \cite{gerdes1988a11}. Figure \ref{mono-linear graphs} shows two mono-linear graphs, that is, they were both made out of only one closed line that embraces all points of the reference frame without repeating itself. 

\subsubsection{Cyclicity}

\begin{figure}
    \centering
    \includegraphics[width=0.65\linewidth]{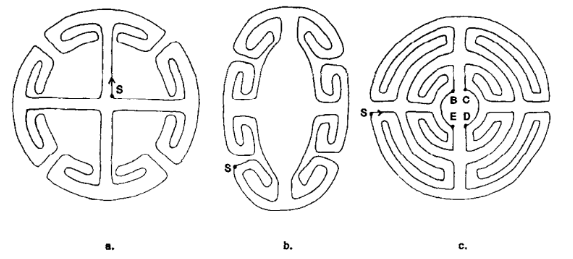}
    \caption{Sand drawing from New Ireland. (a) Four log-drums. (b) The fruit of a wild vine. (c) The thing for which one has been killed. S is the starting and ending point.}
    \label{continuous}
\end{figure}

\begin{figure}
    \centering
    \includegraphics[width=0.25\linewidth]{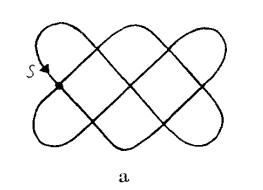}
    \caption{Bushoong sand figures. S denotes the starting of the Bushoong tracings.}
    \label{Bushoong sand figures}
\end{figure}

Upon closer examination of the sand drawings in Malekula, it becomes evident that some are intentionally designed to finish where they begin. There is even a specific term, “suon”, assigned to drawings that conclude at their starting points. Figure \ref{continuous} from Malekula illustrates figures traced continuously, ending where they began \cite{ascher1987a3}. This concept also occurs in the Bushoong sand drawing, as depicted in Figure \ref{Bushoong sand figures} \cite{ascher1988a4}.

\subsubsection{Symmetry}
Symmetry is a prevalent and significant mathematical theme in Western mathematics. Another surprising observation in the corpus is that symmetry is also common in Malekula and among the Tshokwe people. The types of symmetry observed vary, including but not limited to rotational symmetry and bilateral symmetry. Symmetry here is often achieved by the transformation of basic units, such as reflection and inversion.

In Ascher’s study of Malekula sand drawings, symmetry is ubiquitous for even degree graphs. In Table 3 of “Graphs in Cultures: A Study in Ethnomathematics”, there are only 2 out of 52 that lack visual symmetry. Moreover, symmetry with respect to only a single axis is seen in 16 of them. With one exception, the rest (68\%) show double axis symmetry. Four of them show 90° rotational symmetry, and the single exception shows 180° rotational symmetry without double axis symmetry. However, it is opposite for odd degree graphs listed in Table 4 of the paper, where only 3 of 15 are symmetric \cite{ascher1987a3}.

\begin{figure}
    \centering
    \includegraphics[width=0.75\linewidth]{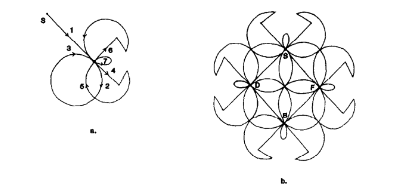}
    \caption{A yam. (a) Procedure A. (b) The result of $AA_{N}A_{R}A_{T}$. (A goes from S to F, $A_{N}$ from F to B, $A_{R}$ from B to D, and $A_{T}$ from D to S.)}
    \label{A yam}
\end{figure}

\begin{figure}
    \centering
    \includegraphics[width=0.5\linewidth]{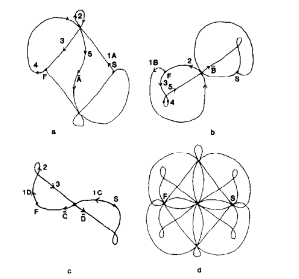}
    \caption{A turtle. (a) Stage 1 with procedure A$\bar{A}$ begins at S and ends at F. (b) Stage 2 with procedure B$\bar{B}$  begins at F and ends at S. (c) Stage 3 with procedure C$\bar{C}$ begins at S and ends at F. Stage 4 with procedure D$\bar{D}$ begins at F and ends at S. (d) The complete path begins and ends at S with procedures A$\bar{A}$B$\bar{B}$C$\bar{C}$D$\bar{D}$.}
    \label{A turtle}
\end{figure}

\begin{figure}
    \centering
    \includegraphics[width=0.5\linewidth]{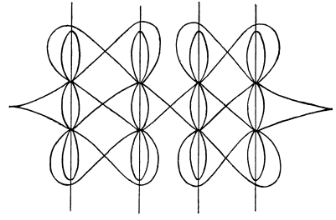}
    \caption{An example of bilateral symmetry.}
    \label{bilateral symmetry}
\end{figure}

Figure \ref{A yam} and Figure \ref{A turtle} show visual rotational symmetries as a result of inversion \cite{ascher1987a3}. Figure \ref{bilateral symmetry} (Hambut hareh navu) \cite{ascher1988a4} is an example of bilateral symmetry where its right and left sides are mirror images of one another. The framework consists of four vertical lines about each one of which another line is drawn so as to form loops; these are then all united by a single line which is interlaced between them, the two ends meeting in a point \cite{deacon1934a6}. 

\begin{figure}
     \centering
     \includegraphics[width=0.5\linewidth]{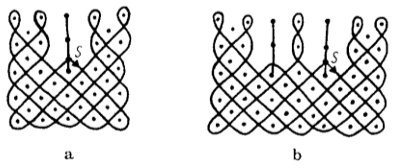}
     \caption{Muyombo trees.}
     \label{Muyombo trees}
 \end{figure}

\begin{figure}
     \centering
     \includegraphics[width=0.5\linewidth]{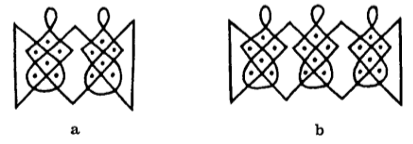}
     \caption{Wild ducks in full flight.}
     \label{Wild ducks}
 \end{figure}

Likewise, in Ascher’s “Graphs in Cultures (II): A Study in Ethnomathematics” which focuses on sand drawings found among the Tshokwe people, the prevalence of symmetry of sona is evident. In Table 1 of the paper, 15 figures have symmetry with respect to only a single axis and 21 have symmetry with respect to a pair of perpendicular axes. Six of the latter also have rotational symmetry of 90°. In addition, 4 have rotational symmetry of 180° including one with rotational symmetry of 90°. Specifically, some symmetrical figures shown in Table 1 are due to the duplication and translation of basic units. Figure \ref{Muyombo trees}b can both be seen as the result of Figure \ref{Muyombo trees}a slit along its line of symmetry and another symmetric strip then inserted. This also applies to Figure \ref{Wild ducks}.

 \begin{figure}
    \centering
    \includegraphics[width=0.5\linewidth]{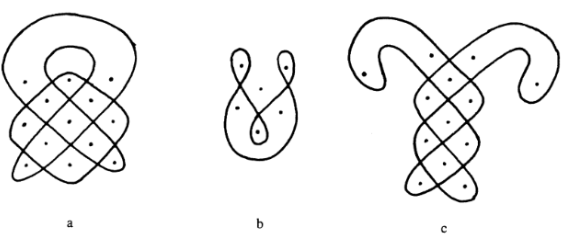}
    \caption{Sona with bilateral symmetries.}
    \label{14 bilateral symmetries}
\end{figure}

\begin{figure}
    \centering
    \includegraphics[width=0.5\linewidth]{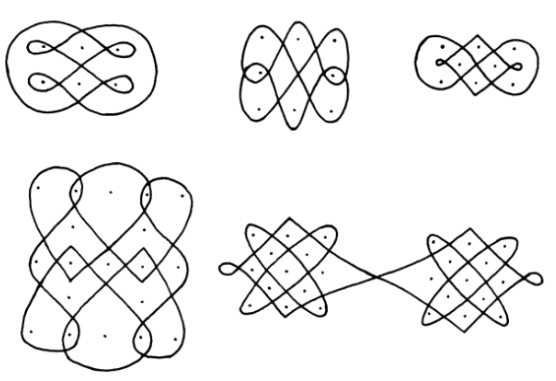}
    \caption{Sona with double bilateral and point symmetries.}
    \label{double bilateral and point symmetries}
\end{figure}

\begin{figure}
    \centering
    \includegraphics[width=0.5\linewidth]{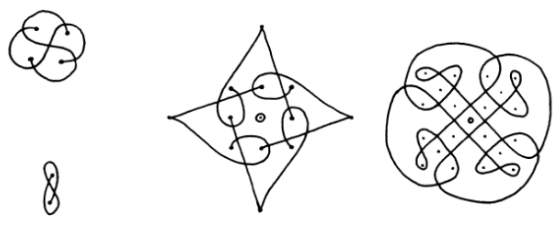}
    \caption{Sona with rotational symmetries.}
    \label{rotational symmetries}
\end{figure}

Paulus Gerdes is another ethnomathematician who had conducted research on sona independently of Ascher. In addition to Ascher who deals with geometrical and topological aspects of sona, in particular with symmetries, extension, enlargement through repetition, and isomorphy, Gerdes concentrates on symmetry and monolinearity (i.e., a whole figure is made up of only one line) \cite{gerdes1994a14}. Gerdes’s work on sona again showcases their bilateral, double bilateral, and rotational symmetry properties. Figure \ref{14 bilateral symmetries} is graphs with bilateral symmetries, while Figure \ref{double bilateral and point symmetries} displays graphs with double bilateral and point symmetry. Additionally, some graphs (e.g. Figure \ref{rotational symmetries}) also exhibit rotational properties, i.e., the corresponding points are at the same distance from the center of rotation \cite{gerdes1988a11}. This all together demonstrates that symmetry is a common sand drawing theme in Malekula and among the Tshokwe people in Angola.

\section{Discussion}

\subsection{Geometrical Attributes}
\begin{figure}
    \centering
    \includegraphics[width=0.6\linewidth]{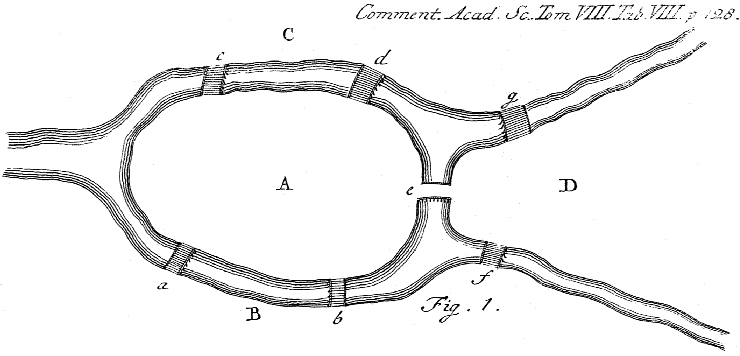}
    \caption{Layout of the seven bridges.}
    \label{Map of Königsberg}
\end{figure}

The continuity concept mentioned in sand drawings perfectly mirrors the idea presented in the notable Königsberg Bridge Problem, the solution of which by Leonhard Euler is considered a foundational moment in the development of graph theory in the Western mathematical world. Essentially, the city of Königsberg in Prussia was situated on both sides of the Pregel River and included two large islands connected to each other and the mainland by seven bridges. The problem was to find a walk through the city that would cross each bridge once and only once (see Figure \ref{Map of Königsberg}). The kind of walk needed to meet the criteria for the problem, i.e. a walk that traverses all edges and repeats none, came to be called the Euler Path \cite{euler1736a17}.

Euler paths are used in the real world when people are planning the best routes to take. For instance, a salesman may wish to find the shortest route through a number of cities and back home again. This gives rise to the travelling salesman problem (TSP), which can be formally stated as follows: Give a finite set of cities N and a distance matrix ($ c_{ij} $) ($i,j \in N$), determine 

\[
\min_{\pi} \left( \sum_{i \in N} c_{i\pi(i)} \right)
\]

where $\pi$ runs over all cyclic permutations of N; $\pi^{k}$(i) is the kth city reached by the salesman from city i. If G denotes the complete directed graph on the vertex set $N$ with a weight $ c_{ij} $ for each arc $(i, j)$, then an optimal tour corresponds to a Hamiltonian circuit on G (i.e. a circuit passing through each vertex exactly once) of minimum total weight \cite{lenstra1975a12}. 

A typical TSP in real life took place in the Institute for Nuclear Physical Research in Amsterdam.  The staff there is working on a challenging issue about the design of computer interfaces. The interfaces consist of modules with multiple pins, and a subset of these pins needs to be interconnected by wires. The positions of the modules are predetermined, and at most two wires are allowed to attach to any pin. The goal is to minimize wire length to prevent signal cross-talk and enhance neatness. Due to the two-wire constraint, the problem is therefore in essence about finding a minimum Hamiltonian path on the graph, which corresponds to a symmetric Euclidean TSP \cite{wilson1996a13}.

The concept of tracing a continuous line through all edges without repetition is evident in New Ireland, Malekula, Bushoong, and Tshokwe cultures. This idea, also a foundational principle in Western graph theory, implies a convergence of mathematical concepts globally. The convergence not only accentuates the need to appreciate and learn from historically undervalued Indigenous cultures but also emphasizes the potential for new discoveries through the exploration of ancient human practices. In this case, we can observe how a simple idea mentioned in the sand drawing relates to the challenging problem of designing high-technology computer interfaces.

\begin{figure}
    \centering
    \includegraphics[width=0.5\linewidth]{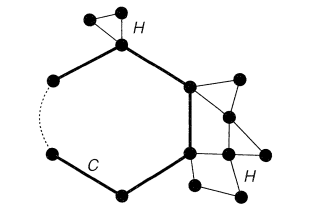}
    \caption{The Euler circuit theorem proof.}
    \label{theorem proof}
\end{figure}

The idea of drawing a graph that has the same starting and ending points relates to another elementary concept in graph theory — Euler circuit. In Euler circuits, each edge is traversed once, and the path forms a closed loop. The Euler Circuit Theorem, the first major result in most graph theory courses, states that a graph possesses an Euler Circuit if and only if the graph is connected and each vertex has even degree. The standard proof may be sketched as follows:

Suppose that P is an Euler path of G. Whenever P passes through a vertex, there is a contribution of 2 towards the degree of that vertex. Since each edge occurs exactly once in P, each vertex must have even degree. The proof is by induction on the number of edges of G. Suppose that the degree of each vertex is even. Since G is connected, each vertex has a degree of at least 2. Therefore, by the Lemma stating that if G is a graph in which the degree of each vertex is at least 2, then G contains a cycle, G contains a cycle C. If C contains every edge of G, the proof is complete. If not, we remove from G the edges of C to form a new, possibly disconnected, graph H with fewer edges than G and in which each vertex still has even degree. By the induction hypothesis, each component of H has an Euler path. Since each component of H has at least one vertex in common with C, by connectedness, we obtain the required Euler path of G by following the edges of C until a non-isolated vertex of H is reached, tracing the Euler path of the component of H that contains that vertex, and then continuing along the edges of C until we reach a vertex belonging to another component of H, and so on. The whole process terminates when we return to the initial vertex (See Figure \ref{theorem proof}) \cite{wilson1996a13}. Again, this can be seen as an example of the unification of important mathematical ideas that occurred in Indigenous culture and the Western mathematical world.

Last but not least, the use of symmetries in Oceanic and African sand drawings reflects the creative imagination of the artists and artisans involved, showcasing their great capacity for abstracting real-life figures. The substantial number of collected graphs containing the symmetry feature suggests that the “drawing experts” who formulated these rules likely understood these rules, i.e. how to make the graphs symmetrically valid. It is suggested that they could even prove some theorems expressed by these rules \cite{meurant1987a15}.

\subsection{Indigenous Roots of Graph Theory}

In the history of mathematics, Euler's solution of the Königsberg bridge problem, published in 1736, is considered to be the first theorem of graph theory and the first true proof in the theory of networks \cite{euler1736a17}. Despite this, the occurrence of analogous mathematical concepts in other Indigenous cultures suggests that intellects from non-European backgrounds may have been aware of and successfully addressed the problem prior to Euler. This raises a similar question to what was proposed by ethnologist Gerhard Kubik \cite{Kubikea30} after he accidentally discovered that one of the sand drawings he collected among the Valuchazi from southern central Africa was identical to a motif on an Indian seal: Is it simply a testimony to the unity of the human mind concerning identical non-numeric mathematical concepts, or does it present a case for diffusion? 

\begin{figure}
    \centering
    \includegraphics[width=0.2\linewidth]{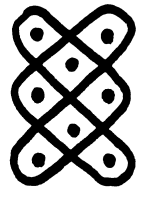}
    \caption{Katuva vufwat}
    \label{Katuva vufwa}
\end{figure}

\begin{figure}
    \centering
    \includegraphics[width=0.5\linewidth]{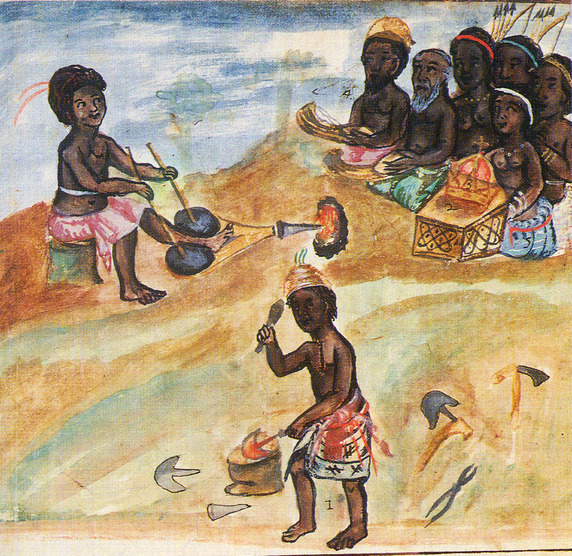}
    \caption{Scene of iron-working in 1650s Angola, with lusona depicted on the upper right, below the crown}
    \label{iron-working}
\end{figure}

By examining how structures and systems of ideas involving number, pattern, logic, and spatial configuration arose in non-European cultures, we can probably explore the origin of graph theory to a further extent. In this case, the fact that sand drawings are produced in sand or clay makes the preservation of sand drawings nearly impossible, resulting in a lack of direct anthropological evidence. Fortunately, there is some circumstantial evidence proving that this tradition must be ancient and certainly pre-colonial in sub-Saharan Africa \cite{Kubikea30}. Early petroglyphs from the Upper Zambezi area in Angola and Citundu-Hulu in the Moçâmedes Desert exhibit structural similarities with lusona ideographs. For example, a lusona known as cingelyengelye and a lusona showing interlaced loops known as zinkhata both appear in the rock arts of the Upper Zambezi recorded by José Redinha. Those petroglyphs date from a period between the 6th century BC and the 1st-century BC \cite{hoddera31}. Therefore, it is possible that those petroglyphs and sona ideographs are related, as evident in the similarities and the geographic location. Additional indirect evidence is found in a watercolor drawing by Italian missionary Antonio Cavazzi de Montecuccolo \cite{Redinhaa32}. Among the drawings in his book is one featuring one of the most fundamental sona, katuva vufwati (see Figure \ref{Katuva vufwa} \cite{Kubika34}), characterized by continuous tracing that concludes where it starts. These illustrations depict scenes of iron-working in 1650s Angola (see Figure \ref{iron-working}). Since Euler published his solution to the Königsberg bridge problem in 1736,  we can infer that the African Indigenous people were aware of what Euler called “problems about the geometry of position” (later becomes graph theory) earlier than Euler. Yet, whether the Indigenous people from Africa succeeded in working out those problems is unknown due to the limited resources in anthropology. 

Nevertheless, even though the Vanuatu islands first had contact with Europeans in April 1606 \cite{borba1990a16}, it was not until 1926, when Deacon carried out ethnographic research in Malekula, that ethnomathematicians started becoming aware of sand drawings in Vanuatu \cite{deacon1934a6}. 

Moreover, it is worthwhile to explore whether Euler drew inspiration from Indigenous practices, such as sand drawings in non-European cultures, to solve the Königsberg bridge problem. Recognizing that the history of mathematics extends beyond European developments is crucial for understanding diverse transmissions of mathematics across cultures, especially when discussing the creation of unified principles like graph theory \cite{Josepha20}. By reviewing Euler’s personal life experiences, we can see that the three important stations of his life were Basel, St. Petersburg, and Berlin — all in Europe. Indeed, Euler received all of his education in Europe, and no evidence suggests he had the chance to have any contact with non-European mathematics \cite{Gautschia21}. Although Euler’s paper on solving the Königsberg bridge problem does not have any references, he mentioned Leibnitz was the first to make mention of the concept of the position geometry. Leibniz, however, treated the East African island as an isolated space for a human laboratory beyond the civilized world \cite{Katza33}. The negative attitude towards African culture implies that Leibniz probably came up with the concept himself without being inspired by any indigenous cultures.

On the other hand, this also stimulates the development of ethnomathematics research, prompting questions like: Has other part of graph theory been contributed by ideas from Indigenous cultures? How have other mathematical areas been influenced by non-European mathematics? 

Another research field related to graph theory worth mentioning is the quipus in Inca culture \cite{Katza28}. Quipus served as a data recording device for Inca leadership, monitoring essential information on taxes, population, required workers for specific projects, etc. A quipu consists of coloured knotted cords, where colours, cord placement, individual cord knots, knot placement, and spaces between knots contribute to the meaning of recorded data. Each quipu has a main cord to which pendant cords are attached, and subsidiary cords may be attached to each pendant cord. Data is recorded on the cords using a system of knots grouped together and separated by spaces, employing a base-10 place-value system. In terms of graph theory, quipus represent a type of graph known as a tree — a connected graph with no cyclic paths, resulting in the number of edges being one less than the number of vertices. Such graphs are often studied in applied algebra courses, as trees are crucial in the modern study of data structures, sorting, and coding theory. In Western mathematics, trees were first defined by Arthur Cayley in 1857. Cayley specifically dealt with rooted trees, where one vertex is designated as the root. He developed a recursive formula to calculate the number of different trees with “r” branches (where “different” is appropriately defined) and applied these results to the study of chemical isomers. 

We lack direct information on whether the Inca quipu makers attempted to answer such questions about their quipus. However, they certainly understood the various ways in which trees could be constructed, as evident in the Incas associating both numbers and colours with each edge of their quipu trees, and the questions they needed to answer in designing them to be useful were not trivial ones. This is similar to what we know about the relationship between the sand drawings discovered in Oceania and Africa and Euler's solution to the Königsberg bridge problem. There is no direct evidence on whether drawers attempted to solve questions like the Königsberg bridge problem, but they were undoubtedly aware of the continuity concept when creating the drawings. Future researchers can possibly conduct a detailed comparative study of how sand drawings correlate with graph theory compared to how quipus relate to graph theory.

\subsection{Educational Potential}
But after examining the relationship between sand drawings in Oceania and Africa and the graph theory developed in Europe, we might ask ourselves: What exactly is mathematics? 

Within the Western education system, mathematics is often described as people's least favorite subject, and, in fact, people sometimes say, without embarrassment, that they are “no good” at mathematics \cite{Latterella23}. In a study involving 74 future elementary education teachers enrolled in a “Mathematics for Elementary Education Majors” course at the University of Minnesota Duluth in 2009, students were asked to write down what they thought mathematics was. Surprisingly, the vast majority of responses were that mathematics is about learning to perform mathematical operations (adding, subtracting, etc.). A representative response was, “Mathematics is the adding, subtracting, etc. of numbers.” Undoubtedly, this oversimplified and narrow conception of mathematics could influence their future students' perception of the subject \cite{Latterella22}. Hence, we can observe that despite the Western mathematics education system having cultivated numerous brilliant mathematicians, there are still improvements that educators can make to enable students to better appreciate the beauty of mathematics and be aware of the strong correlation between mathematics and our daily practices.

It is, therefore, time for us to shift our focus to mathematics from other cultures to see what we might learn from them. The Ministry of Education and the Vanuatu Cultural Center prompt Ni-Vanuatu as well as foreign researchers to take part in the elaboration of teaching programs related to social and cultural concerns, through their participation in conferences/workshops/trainings organized in the capital Port Vila in particular. Regarding the teaching of mathematics, Georges Tauanearu, Professor of mathematics at the “Vanuatu Institute of Teacher Education” — VITE teacher training, has been giving a course in ethnomathematics for about a decade. This course allows students to study seminal works in ethnomathematics, while being encouraged to use as pedagogical materials — the local mathematical practices that they would be able to collate in the cultural/linguistic areas where they will be sent to teach in local secondary schools, throughout the archipelago.

\begin{figure}
    \centering
    \includegraphics[width=1\linewidth]{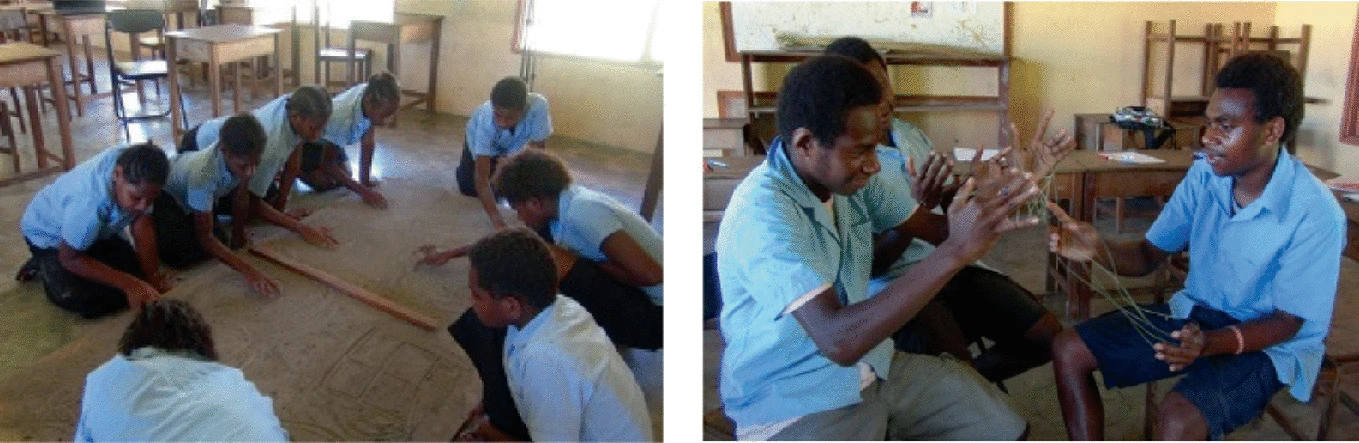}
    \caption{Bringing string figure making and sand drawing into the classroom, Topol High School, North Ambrym, Vanuatu.
}
    \label{classroom}
\end{figure}

In 2019, ethnomathematician Eric Vandendriessche \cite{Vandendriessche27} collaborated with a mathematics teacher Frederick Worwor from Vanuatu to conduct a sequence of three teaching sessions in the 6th and 10th-year classrooms on string figure-making and sand drawing  (See Figure \ref{classroom}). For the making of sand drawings, a couple of wooden trays and some ashes spread on them had been prepared for the students. The activities involved encouraging students to produce their own sand drawings and leading students to identify and analyze mathematical concepts behind the drawings, such as iteration and alteration. The homework assigned focused on sharing acquired procedures and identifying local experts in villages. Overall, the experiment can be considered a good attempt to integrate local activities with mathematical characteristics into the curriculum.

\begin{figure}
    \centering
    \includegraphics[width=0.2\linewidth]{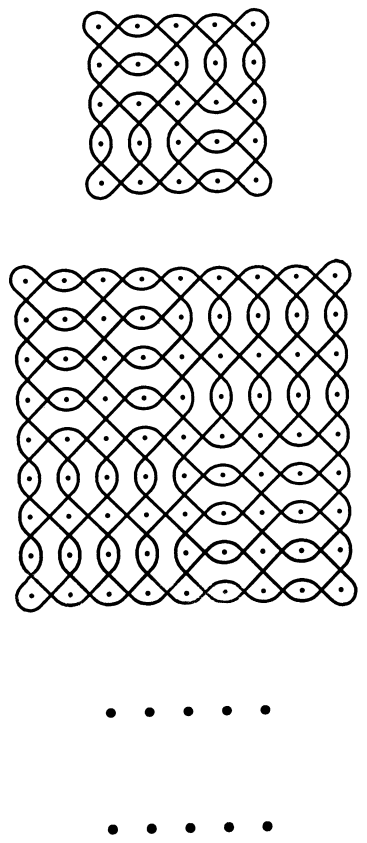}
    \caption{“Find the missing figures” problem.}
    \label{missing figures}
\end{figure}

On the other hand, the potential of using the Tshokwe sand drawings in mathematical education has been explored by Gerdes in his paper “On possible uses of traditional Angolan sand drawings in the mathematics classroom” published in 1988 \cite{gerdes1988a11}. The examples include the study of arithmetical relationships and progressions, symmetry, similarity, Euler graphs, and greatest common divisors using sand drawings. As a variant on the well-known theme of arithmetical problems of the type “Find the missing numbers”, a series of geometric problems has been elaborated: “Find the missing figures”. Figure \ref{missing figures} provides an example. Given two sona of a series, students are asked to discover the third and subsequent sona. To find the solutions, both the dimensions of the successive drawings and the geometrical algorithm have to be carefully observed.  These problems aim to develop a sense of geometric algorithms, generalization, and symmetry \cite{gerdes1988a11}. 

Future researchers can work on similar methodologies to integrate diverse forms of ethnomathematics into the curriculum, aiming to not only apply ethnomathematics study findings to real-world scenarios but also enrich students' understanding of the subject. Furthermore, for children from ethnic minority communities, their cultural identity can be better recognized, and the multicultural nature of these societies in the school curriculum can be more accurately reflected. 

\begin{figure} 
    \centering
    \includegraphics[width=0.8\linewidth]{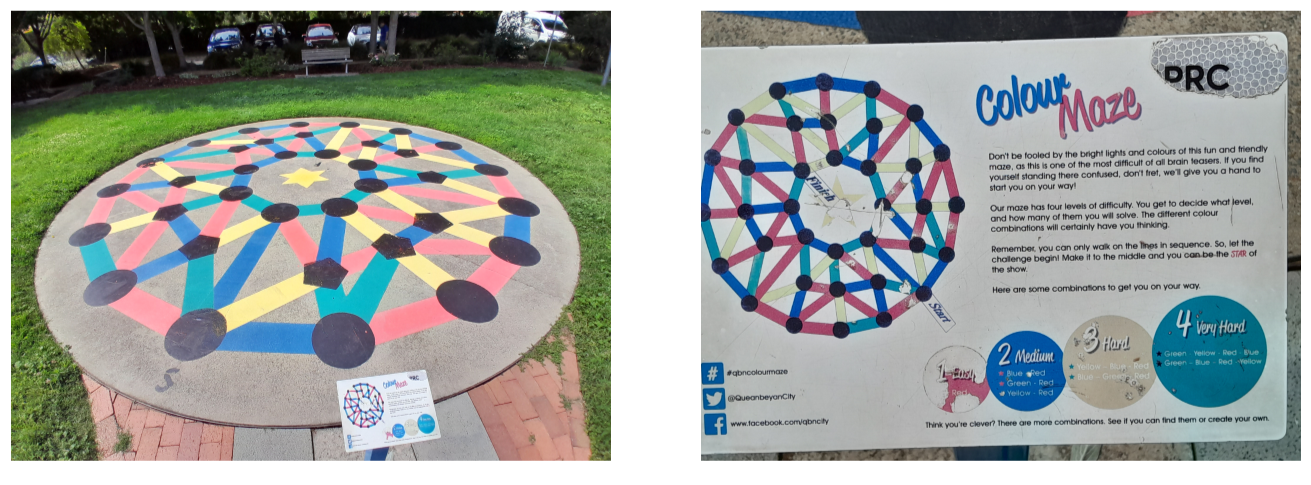}
    \caption{A colour maze in Queanbeyan park.}
    \label{colour maze}
\end{figure}

Having said that, the implementation of sand drawings in the classroom to support students' learning processes is a deliberate strategy employed by educators to bring attention to real-life non-numeric mathematics. It is essential for readers themselves to recognize that mathematics is closely intertwined with our daily lives, and often, we simply lack the awareness to fully embrace the diverse ways in which mathematics manifests across various cultures. Below are two photos (Figure \ref{colour maze}) depicting a colour maze in a park in Queanbeyan, a city in the south-eastern region of New South Wales, Australia. These images serve as a classic example of graph theory in real life, encouraging readers to explore the role of mathematics in their own lives, understanding that mathematics goes beyond algebraic calculations from the textbook.

\section{Conclusion}

An introduction to ethnomathematics and the history that led to Ascher’s studies on sand drawings in Oceania islands and Africa has been presented at the start of the paper. Ethnomathematics, as defined by Ubiratan D'Ambrosio in 1985 \cite{borba1990a16}, is the mathematics practiced among identifiable cultural groups, such as national-tribal societies, labour groups, children of a certain age bracket, professional classes, and so on. Even the mathematics produced by professional mathematicians can be seen as a form of ethnomathematics because it was produced by an identifiable cultural group and because it is not the only mathematics that has been produced. Based on this statement, we can infer that certain mathematical ideas found in traditional practices, such as sand drawings, hold equal significance to those in graph theory taught in schools. This raises our awareness of how we should treat sand drawings with more care and reflect on the possible geometric or topological ideas linked behind them. Even though the resemblance is not exact, the fundamental concepts of the related graph theory have been highlighted in the paper so that readers who are interested in further exploring this topic may have a quick reference.

Three illustrative examples have been presented in the Results and Discussion section. The concept of tracing figures continuously in the sand is closely associated with Euler’s path and the Traveling Salesman Problem. Additionally, the goal of drawing a graph that ends where it begins represents an essential aspect of the Euler circuit. Lastly, the prevalence of diverse forms of symmetry is meticulously examined and compared between sand drawings found in Malekula and those among the Tshokwe people. The discussion briefly touches upon the Indigenous origins of graph theory and the educational potential found in sand drawings. Future researchers are encouraged to explore these fields more extensively to contribute to their exploration and understanding. Motivated by Ascher’s work on sand drawings, this paper provides solid examples of how to relate some mathematical ideas contained in aboriginal practices to modern mathematics. It is hoped that this can capture readers' attention to our day-to-day mathematical experiences and realise that mathematics is an integral part of the shared human experience. 

\section{Acknowledgment}
This paper would not have been possible without the guidance and help of my two supervisors, Rowena Ball and Hongzhang Xu, who contributed their valuable assistance in the preparation and completion of this study. I extend my sincere gratitude to Rowena Ball for generously providing the captivating colour maze photo in Queanbeyan Park, which significantly enriched the visual aspects of this paper. Special thanks are also due to Hongzhang Xu for his invaluable feedback and meticulous comments on each version of the paper. His insights and attention to detail have significantly improved the overall quality and clarity of the work.

\printbibliography

\clearpage
\thispagestyle{empty}
\listoffigures

\end{document}